\documentclass{amsart}[11pt,amssymb]

\usepackage{amsmath,amssymb,amsthm,amsfonts,hyperref}

\DeclareMathSymbol{\twoheadrightarrow} {\mathrel}{AMSa}{"10}

\def\Q{{\mathbf Q}}

\def\Z{{\mathbf Z}}

\def\F{{\mathbf F}}

\def\Gal{\mathrm{Gal}}

\def\End{\mathrm{End}}

\def\Aut{\mathrm{Aut}}

\def\Hom{\mathrm{Hom}}

\def\I{\mathrm{Id}}

\def\fchar{\mathrm{char}}

                                    \def\tr{\mathrm{tr}}

\def\Gp{\mathrm{Gp}}

                \def\sL{\mathfrak{sl}}

        \def\K_a{\bar{K}}

\def\dim{\mathrm{dim}}

\def\g{{\mathfrak g}}

\def\K{{\mathcal{K}}}

\newtheorem{thm}{Theorem}[section]

\newtheorem{lem}[thm]{Lemma}

\newtheorem{cor}[thm]{Corollary}

\newtheorem{prop}[thm]{Proposition}

\theoremstyle{definition}

\newtheorem{rem}[thm]{Remark}

\title[Endomorphisms of abelian varieties]{Endomorphisms of abelian varieties,  cyclotomic extensions and Lie algebras}

\author{Yuri G. Zarhin}

\address{Department of Mathematics, Pennsylvania
State University, University Park, PA 16802, USA}

\address{Institute for Mathematical Problems in Biology,
Russian Academy of Sciences, Pushchino, Moscow Region, Russia}

\email{zarhin\char`\@math.psu.edu}

\begin{document}
\begin{abstract}
We prove an analogue of the Tate conjecture on homomorphisms of
abelian varieties over infinite cyclotomic extensions of finitely
generated fields of characteristic zero.
\end{abstract}

\maketitle

\section{Introduction}

The aim of this note is to extend Faltings' results \cite{F1,F2}
concerning the Tate conjecture on homomorphisms of abelian varieties
\cite{Tate1,Tate2} over finitely generated fields $K$ of
characteristic zero to their infinite cyclotomic extensions
$K(\ell)=K(\mu_{\ell^{\infty}})$. The possibility of such
generarization  (in the case of number fields $K$) was stated
(without a detailed proof) in \cite[\S 6, Subsect. {\bf 2}, pp.
91--92]{ZarhinParshin}; it was pointed out there that this result
follows from the theorem of Faltings combined with technique
developed in \cite{ZarhinTorsion77} and a theorem of F.A. Bogomolov
about homotheties \cite{Bogomolov1,Bogomolov2}. Our main result is
the following assertion. (Here $\Gal(E)$ stands for the absolute
Galois group of $E$ while $T_{\ell}(X)$ and $T_{\ell}(Y)$ are the
Tate modules of abelian varieties of $X$ and $Y$ respectively.)

\begin{thm}
\label{main} Suppose that $K$ is a field that is finitely generated over  $\Q$
and $\ell$ is a prime. Let us put $E=K(\ell)$. If $X$ and $Y$ are abelian
varieties over $E$ then the natural embedding of $\Z_{\ell}$-modules
$$\Hom_E(X,Y)\otimes \Z_{\ell} \hookrightarrow \Hom_{\Gal(E)}(T_{\ell}(X),
T_{\ell}(Y))$$ is bijective.
\end{thm}

\begin{rem}
\label{enlarge} Replacing $K$ by its suitable finite (sub) extension
(of $E$), we may in the course of the proof of Theorem \ref{main}
assume that both $X$ and $Y$ are defined over $K$.
\end{rem}

\begin{rem}

A.N. Parshin \cite[\S 6, Subsect. {\bf 2}, pp.
91--92]{ZarhinParshin} conjectured that the following analogue of
the Mordell conjecture holds true: if $K$ is a number field and $C$
is an absolutely irreducible smooth projective curve over
$E=K(\ell)$ then the set $C(E)$ of its $E$-rational points is finite
if the genus of $C$ is greater than $1$. Theorem \ref{main} has
arisen from attempts to understand which parts of Faltings' proof of
the Mordell conjecture \cite{F1} can be extended to the case of
$K(\ell)$.
\end{rem}

The paper is organized as follows. In Section \ref{Intro} we discuss the
$\ell$-adic Lie algebras arising from Tate modules of abelian varieties and
their centralizers. In Section \ref{homo} we deal with analogues of the Tate
conjecture on homomorphisms over arbitrary fields. In Section \ref{central} we
prove the main result.


I am grateful to A.N. Parshin for stimulating discussions.

\section{Tate modules, $\ell$-adic Lie groups and Lie algebras}
\label{Intro}

Let $K$ be a field, $\bar{K}$ its algebraic closure and
$\Gal(K)=\Aut(\bar{K}/K)$ its absolute Galois group. If $m$ is a
positive integer that is not divisible by $\fchar(K)$ then we write
$\mu_m$ for the cyclic order $m$ multiplicative subgroup of $m$th
roots of unity in $\bar{K}$ and $K(\mu_m)$ for the corresponding
cyclotomic extension of $K$. If $\ell$ is a prime different from
$\fchar(K)$ then we write $K(\ell)$ for the ``infinite" cyclotomic
extension
$$E=K(\mu_{\ell^{\infty}})=\cup_{i=1}^{\infty}K(\mu_{\ell^i}).$$
 It is
well known that the compact  Galois group $\Gal(K(\ell)/K)$ is
canonically isomorphic to a closed subgroup of $\Z_{\ell}^*$. We
write $$\chi_{\ell}:\Gal(K) \to \Z_{\ell}^*$$ for the corresponding
{\sl cyclotomic} character; its kernel coincides with
$\Gal(K(\ell))=\Aut(\bar{K}/K(\ell))$.

We write $\Z_{\ell}(1)$ for the projective limit of the groups
$\mu_{\ell^j}$ where the transition map $\mu_{\ell^{j+1}}\to
\mu_{\ell^j}$ is raising to $\ell$th power. It is well known that
 $\Z_{\ell}(1)$ is a free $\Z_{\ell}$-module of rank
 $1$ provided with the natural structure of a Galois module while
 the defining homomorphism
$$\Gal(K) \to \Aut_{\Z_{\ell}}(\Z_{\ell}(1))=\Z_{\ell}^{*}$$
coincides with  $\chi_{\ell}$. Let us consider the $1$-dimensional
$\Q_{\ell}$-vector space
$$\Q_{\ell}(1)=\Z_{\ell}(1)\otimes_{\Z_{\ell}} {\Q_{\ell}},$$
which carries the structure of a Galois module provided by the same
character
$$\chi_{\ell}:\Gal(K) \to  \Z_{\ell}^{*}\subset
\Q_{\ell}^{*}=\Aut_{\Q_{\ell}}(\Q_{\ell}(1)).$$

Let $A$ be an abelian variety over $K$ and $\End_K(A)$ the ring of its
$K$-endomorphisms. If $X$ and $Y$ are abelian varieties over $K$ then we write
$\Hom_K(X,Y)$ for the group of $K$-homomorphisms from $X$ to $Y$. If $m$ is as
above  then we write $A_m$ for the kernel of multiplication by $m$ in
$A(\bar{K})$. The subgroup $A_m$ is a free $\Z/m\Z$-module of rank $2\dim(A)$
\cite{MumfordAV} and a Galois submodule of $A(\bar{K})$.

We write $T_{\ell}(A)$ for the $\ell$-adic Tate module of $A$, which is the
projective limit  of the groups $A_{\ell^j}$ while the transition map
$A_{\ell^{j+1}}\to A_{\ell^j}$ is multiplication by $\ell$ \cite{MumfordAV}. It
is well-known that $T_{\ell}(A)$ is a free $\Z_{\ell}$-module of rank
$2\dim(A)$, the natural map $T_{\ell}(A) \to A_{\ell^j}$ gives rise to the
isomorphisms
$$T_{\ell}(A)/\ell^j =A_{\ell^j}$$
and the Galois actions on $A_{\ell^j}$  give rise to
  the natural continuous homomorphism ($\ell$-{\sl adic representation})
$$\rho_{\ell,A}=\rho_{\ell,A,K}: \Gal(K) \to \Aut_{\Z_{\ell}}(T_{\ell}(A)),$$
which provides $T_{\ell}(A)$ with the natural structure of a
$\Gal(K)$-module \cite{SerreAbelian}. The image
$$G_{\ell,A,K} =\rho_{\ell,A}(\Gal(K)) \subset \Aut_{\Z_{\ell}}(T_{\ell}(A))$$
is a compact $\ell$-adic Lie group \cite{SerreAbelian}; clearly,
$G_{\ell,A,K}\subset 1+\ell^j \End_{\Z_{\ell}}(T_{\ell}(A))$ if and
only if $\Gal(K)$ acts trivially on $T_{\ell}(A)/\ell^j
=A_{\ell^j}$, i.e., $A_{\ell^j}\subset A({K})$.

Clearly, $\Gal(K(\ell))=\Aut(\bar{K}/K(\ell))$ is a compact normal
subgroup of $\Gal(K)$. We write $G_{\ell,A,K}^{1}$ for its image
$\rho_{\ell,X}(\Gal(K(\ell)))$, which is a compact normal Lie
subgroup of $G_{\ell,A,K}$. By definition,
$$G_{\ell,A,K}^{1}=G_{\ell,A,K(\ell)}.$$

Let us consider the $2\dim(A)$-dimensional $\Q_{\ell}$-vector space
$$V_{\ell}(A)=T_{\ell}(A)\otimes_{\Z_{\ell}}\Q_{\ell}.$$
One may view $T_{\ell}(A)$ as the $\Z_{\ell}$-lattice in $V_{\ell}(A)$. We have
$$G_{\ell,A,K}\subset \Aut_{\Z_{\ell}}(T_{\ell}(A))\subset \Aut_{\Q_{\ell}}(V_{\ell}(A)),$$
which provides $V_{\ell}(A)$ with the natural structure of a $\Gal(K)$-module.

Notice that the Lie algebra of the compact $\ell$-adic Lie group
$\Aut_{\Z_{\ell}}(T_{\ell}(A))$ coincides with $\End_{\Q_{\ell}}(V_{\ell}(A))$.
The Lie algebra $\g_{\ell,A}$ of $G_{\ell,A,K}$ is a $\Q_{\ell}$-linear Lie
subalgebra of $\End_{\Q_{\ell}}(V_{\ell}(A))$. The Lie algebra $\g_{\ell,A}^0$
of $G_{\ell,A,K}^{1}$ is an ideal in $\g_{\ell,A}$. It is known
\cite{SerreAbelian} that the Lie algebras $\g_{\ell,A}$ and $\g_{\ell,A}^0$
will not change if we replace $K$ by its finite algebraic extension.

Let $\I$ be the identity map on $V_{\ell}(A)$ and let
$$\tr: \End_{\Q_{\ell}}(V_{\ell}(A)) \to \Q_{\ell}$$ be the trace map. Let
$$\det:\Aut_{\Q_{\ell}}(V_{\ell}(A)) \to \Q_{\ell}^{*}$$ be the determinant map.
We write $\sL(V_{\ell}(A))$ for the Lie subalgebra of traceless
linear operators in $\End_{\Q_{\ell}}(V_{\ell}(A))$.

Let $\End_{\g_{\ell,A}} (V_{\ell}(A))$ be the the centralizer of $\g_{\ell,A}$
in $\End_{\Q_{\ell}}(V_{\ell}(A))$ and $\End_{\g_{\ell,A}^0} (V_{\ell}(A))$ be
the the centralizer of $\g_{\ell,A}^0$ in $\End_{\Q_{\ell}}(V_{\ell}(A))$.
Clearly,
$$\Q_{\ell}\I\subset \End_{\g_{\ell,A}} (V_{\ell}(A))\subset \End_{\g_{\ell,A}^0}
(V_{\ell}(A))\subset \End_{\Q_{\ell}}(V_{\ell}(A)).$$

\begin{rem}
\label{Liecentral} Since $\g_{\ell,A}$ is the Lie algebra of
$G_{\ell,A,K}$,
$$\End_{\g_{\ell,A}} (V_{\ell}(A)) \supset \End_{G_{\ell,A,K}}(V_{\ell}(A))=\End_{\Gal(K)}(V_{\ell}(A)).$$
Since $\g_{\ell,A}^0$ is
the Lie algebra of $G_{\ell,A,K}^{1}$,
$$\End_{\g_{\ell,A}^0} (V_{\ell}(A)) \supset \End_{G_{\ell,A,K}^1}(V_{\ell}(A))=
\End_{\Gal(K(\ell))}(V_{\ell}(A))=\End_{\Gal(E)}(V_{\ell}(A)).$$
\end{rem}

In  the following two propositions we assume that $A$ has positive
dimension.

\begin{prop}
\label{centralizerLie}
\begin{enumerate}
\item $\g_{\ell,A}^0=\g_{\ell,A} \cap \sL(V_{\ell}(A))\subset
\sL(V_{\ell}(A))$.

\item Suppose that $\g_{\ell,A}$ contains the homotheties $\Q_{\ell}\I$. Then
$$\g_{\ell,A}=\Q_{\ell}\I\oplus \g_{\ell,A}^0.$$ In particular, the centralizers
$\End_{\g_{\ell,A}} (V_{\ell}(A))$  and  $\End_{\g_{\ell,A}^0} (V_{\ell}(A))$
do coincide.
\end{enumerate}
\end{prop}

\begin{proof}
It is known \cite[Sect. 1.3]{SerreIZV} that a choice of a polarization on $A$
gives rise to a nondegenerate alternating bilinear form
$$e_l: V_{\ell}(A) \times V_{\ell}(A) \to \Q_{\ell}(1)$$
such that
$$e_l(\rho_{\ell,A}(\sigma)(x), \rho_{\ell,A}(\sigma)(y))=\chi_{\ell}(\sigma)\cdot e_{\ell}(x,y) \
\forall x,y \in V_{\ell}(A); \sigma \in \Gal(K).$$ By fixing a
(non-canonical) isomorphism of $\Q_{\ell}$-vector spaces
$$\Q_{\ell}(1) \cong \Q_{\ell},$$
we may assume that the alternating form $e_l$ takes on values in
$\Q_{\ell}$. We obtain that $G_{\ell,A,K}$ lies in the corresponding
group of symplectic similitudes
\begin{multline*}
$$\Gp(V_{\ell}(A),e_{\ell}) =\{s \in \Aut_{\Q_{\ell}}(V_{\ell}(A))\mid \exists c \in \Q_{\ell}^{*} \ \mbox{such that }\\
 e_l(sx,sy)=
c\cdot  e_l(x,y) \ \forall x,y\in V_{\ell}(A)\}\subset
\Aut_{\Q_{\ell}}(V_{\ell}(A))
\end{multline*}
 and $\chi_{\ell}$ coincides with the
composition of
$$\rho_{\ell,A}:\Gal(K) \twoheadrightarrow G_{\ell,A,K}$$ and
$$G_{\ell,A,K}\subset \Gp(V_{\ell}(A),e_{\ell}) \stackrel{\mathrm{c}}{\to} \Q_{\ell}^{*}$$
where the scalar $\mathrm{c}(s)$ is defined by
$$e_{\ell}(sx,sy) =\mathrm{c}(s)\cdot e_{\ell}(x,y) \ \forall x,y\in  V_{\ell}(A); s \in \Gp(V_{\ell}(A),e_{\ell}).$$
Clearly, $$\mathrm{c}:\Gp(V_{\ell}(A),e_{\ell}) {\to} \Q_{\ell}^{*},
\ s \mapsto \mathrm{c}(s)$$ is a homomorphism of $\ell$-adic Lie
groups and
$$\mathrm{c}(s)^{\dim(A)}=\det(s) \ \forall s \in \Gp(V_{\ell}(A),e_{\ell})$$
(recall that $V_{\ell}(A)$ is a $2\dim(A)$-dimensional
$\Q_{\ell}$-vector space). It is also clear that $G_{\ell,A,K}^{0}$
coincides with the kernel of  the homomorphism of $\ell$-adic Lie
groups
$$\mathrm{c}: G_{\ell,A,K}\subset \Gp(V_{\ell}(A),e_{\ell})\to \Q_{\ell}^{*},$$
and therefore $\g_{\ell,A}^0$ coincides with the kernel of the
corresponding tangent map of Lie algebras
$$\g_{\ell,A} \to \Q_{\ell}.$$
On the other hand, one may easily check that the tangent map is
$$\frac{1}{\dim(A)}\tr: \g_{\ell,A} \to \Q_{\ell}$$
(because $\tr$ is the tangent map to $\det$.) This implies that
$$\g_{\ell,A}^0=\g_{\ell,A} \cap \sL(V_{\ell}(A)).$$
This proves the first assertion of Proposition. In order to prove the second
assertion, notice that $\Q_{\ell}\I \cap \sL(V_{\ell}(A))=\{0\}$ and therefore
$\g_{\ell,A}$ contains $\Q_{\ell}\I\oplus\g_{\ell,A}^0$. On the other hand,
since $\g_{\ell,A}^0$ is the kernel of $\g_{\ell,A} \to \Q_{\ell}$, its
codimension in $\g_{\ell,A}$ is (at most) $1$. This implies that $\g_{\ell,A}=
\Q_{\ell}\I\oplus\g_{\ell,A}^0$.
\end{proof}

\begin{prop}
\label{centralizerGroup} There exists a finite separable algebraic field
extension $K_0/K$ that enjoys the following properties.

If $K^{\prime}/K_0$ is a  finite separable algebraic field extension
then
$$\End_{\Gal(K^{\prime})}(V_{\ell}(A))=\End_{\g_{\ell,A}} (V_{\ell}(A)).$$
\end{prop}

\begin{proof} (Compare with \cite[Prop. 1 and its proof]{SerreLF}.)
Let us choose  open neighborhoods $V$ of $0$ in $\g_{\ell,A}$ and
$U$ of $\I$ in $G_{\ell,A,K}$ such that the $\ell$-adic exponential
map $\exp$
 and  logarithm map $\log$ establish mutually inverse $\Q_{\ell}$-analytic
isomorphisms between $V$ and $U$.

Let $G_0$ be an open subgroup of $G_{\ell,A,K}$ that lies in $U$.
(The existence of such $G_0$ follows from Corollary 2 in \cite[Part
II, Ch. 4, Sect. 8, p. 117]{SerreLie}.)

 Then $V_0=\log(G_0)$ is an open subset of $\g_{\ell,A}$ that
contains $0$ and $G_0=\exp(V_0)$.
 Clearly, $G_0$ has finite index in $G_{\ell,A,K}$ and
$$\End_{G_0}(V_{\ell}(A))=\End_{V_0}(V_{\ell}(A))=\End_{\g_{\ell,A}} (V_{\ell}(A)).$$

 The preimage of $G_0$ in $\Gal(K)$ is an open subgroup of finite index and
therefore coincides with $\Gal(K_0)$ for a certain finite separable algebraic
field extension $K_0$ of $K$. It follows that
$$\End_{\Gal(K_0)}(V_{\ell}(A))=\End_{G_0}(V_{\ell}(A))=\End_{\g_{\ell,A}} (V_{\ell}(A)).$$
If $K^{\prime}/K_0$ is a  finite separable algebraic field extension
then $\Gal(K^{\prime})$ is a compact subgroup of finite index in
$\Gal(K_0)$ and its image
$G^{\prime}=\rho_{\ell,A}(\Gal(K^{\prime}))$ is a closed subgroup of
finite index in $G_0$ and therefore is open in $G_0$ and therefore
is also open in $G_{\ell,A,K}$.
As above, $V^{\prime}=\log(G^{\prime})$ is an open subset of
$\g_{\ell,A}$ that contains $0$ and $G^{\prime}=\exp(V^{\prime})$
and
$$\End_{\Gal(K^{\prime})}(V_{\ell}(A))=\End_{G^{\prime}}(V_{\ell}(A))=\End_{V^{\prime}}(V_{\ell}(A))=\End_{\g_{\ell,A}} (V_{\ell}(A)).$$
\end{proof}

\section{Homomorphisms of abelian varieties}
\label{homo} Throughout this Section, $X,Y,Z,A$ are abelian
varieties over $K$ and $\ell$ is a prime different from $\fchar(K)$.
We write $\Hom_K(X,Y)$ for the (finitely generated free) commutative
group of $K$-homomorphisms from $X$ to $Y$. If $X=Y=A$ then
$\Hom_K(X,Y)$ coincides with the ring $\End_K(A)$ of
$K$-endomorphisms of $A$.

There is a natural embedding of $\Z_{\ell}$-modules
$$i_{X,Y,K}:\Hom_K(X,Y)\otimes \Z_{\ell} \hookrightarrow \Hom_{\Gal(K)}(T_{\ell}(X),
T_{\ell}(Y)),$$ whose cokernel is torsion free \cite{MumfordAV,Tate2}.

\begin{rem}
\label{extn} Notice that if $L/K$ is a finite or infinite Galois
extension and
$$i_{X,Y,L}:\Hom_L(X,Y)\otimes \Z_{\ell} \hookrightarrow \Hom_{\Gal(L)}(T_{\ell}(X),
T_{\ell}(Y))$$ is bijective then $i_{X,Y,K}$ is also bijective. This assertion
follows easily from the following obvious description of $\Gal(L/K)$-invariants
$$\Hom_K(X,Y)=\{\Hom_L(X,Y)\}^{\Gal(L/K)},$$
$$\Hom_{\Gal(K)}(T_{\ell}(X), T_{\ell}(Y))=\{\Hom_{\Gal(L)}(T_{\ell}(X),
T_{\ell}(Y))\}^{\Gal(L/K)}$$ and the $\Gal(L/K)$-equivariance of $i_{X,Y,L}$.
\end{rem}

Extending $i_{X,Y,K}$ by $\Q_{\ell}$-linearity, we obtain the natural embedding
of $\Q_{\ell}$-vector spaces
$$\tilde{i}_{X,Y,K}:\Hom_K(X,Y)\otimes \Q_{\ell} \hookrightarrow \Hom_{\Gal(K)}(V_{\ell}(X),
V_{\ell}(Y)),$$ see \cite[Section 1, displayed formula (2) on p. 135]{Tate2}.

The following observations are due to  J. Tate \cite[Sect. 1, Lemma 1 and Lemma
3 and its proof on p. 135]{Tate2}.
\begin{lem}[of Tate]
\label{prod}
\begin{enumerate}
\item The map $i_{X,Y,K}$ is bijective if and only $\tilde{i}_{X,Y,K}$ is
bijective.

\item

Let us put $Z=X\times Y$. If the embedding
$$\tilde{i}_{Z,Z,K}:\End_K(Z)\otimes\Q_{\ell} \hookrightarrow
\End_{\Gal(K)}(V_{\ell}(Z))$$ is bijective then
$$\tilde{i}_{X,Y,K}:\Hom_K(X,Y)\otimes \Q_{\ell} \hookrightarrow \Hom_{\Gal(K)}(V_{\ell}(X),
V_{\ell}(Y))$$ is also bijective.
\end{enumerate}
\end{lem}

Tate \cite{Tate1,Tate2} conjectured and G. Faltings \cite{F1,F2} proved that
this embedding is actually a bijection when $K$ is finitely generated over the
field $\Q$ of rational numbers.

\begin{thm}
\label{main0} Suppose that $K$ is field  and $\ell$ is a prime that is
different from $\fchar(K)$. Let us put $E=K(\ell)$. Suppose that $A$ is an
abelian variety of positive dimension over $K$ such that for all finite
separable algebraic field extensions $K^{\prime}/K$ the embedding
$$\tilde{i}_{A,A,K^{\prime}}:\End_{K^{\prime}}(A)\otimes\Q_{\ell} \hookrightarrow
\End_{\Gal(K^{\prime})}(V_{\ell}(A))$$ is bijective. If
$\g_{\ell,A}$ contains the homotheties $\Q_{\ell}\I$ then the
injective maps
$$\tilde{i}_{A,A,E}:\End_{E}(A)\otimes\Q_{\ell} \hookrightarrow
\End_{\Gal(E)}(V_{\ell}(A))$$ and
$$i_{A,A,E}:\End_{E}(A)\otimes \Z_{\ell} \hookrightarrow
\End_{\Gal(E)}(T_{\ell}(A))$$ are bijective.
\end{thm}

\begin{proof}
Let us consider the field  $K_2=K(A_{\ell^2})$ of definition of all
points of $A_{\ell^2}$ and put  $E_2=K_2(\ell)$. Clearly, $K_2/K$
and $E_2/E$ are finite Galois extensions. Applying Remark \ref{extn}
(to $E_2/E$ instead of $L/K$) and Lemma \ref{prod}.1 to $X=Y=A$, we
observe that in the course of the proof we may (and will) assume
that
$$K=K_2=K(A_{\ell^2}),$$
i.e., $A_{\ell^2}\subset A(K)$.
Since $\ell^2 \ge 4$, it follows from a result of A. Silverberg
\cite{Silverberg} that all $\bar{K}$-endomorphisms of $A$ are
defined over $K$. In particular,
$$\End_K(A)=\End_E(A)=\End_{\bar{K}}(A).$$ Using Proposition
\ref{centralizerGroup}, we may replace $K$ by its finite separable
algebraic extension in such a way that
$$\End_{\Gal(K^{\prime})}(V_{\ell}(A))=\End_{\g_{\ell,A}} (V_{\ell}(A))$$ for
all finite separable algebraic field extensions $K^{\prime}$ of $K$.

Let $K_0/K$ satisfies the conclusion of Proposition
\ref{centralizerGroup}. Replacing $K_0$ by its normal closure over
$K$, we may and will assume that $K_0/K$ is a finite Galois
extension. Let us put $E_0=K_0(\ell)$. Clearly, $E_0/E$ is a finite
Galois extension and
$$\End_E(A)=\End_K(A)=\End_{K_0}(A)=\End_{E_0}(A).$$
By the assumption of Theorem \ref{main0},
$$\End_{\Gal(K_0)}(V_{\ell}(A))=\End_{K_0}(A)\otimes\Q_{\ell}.$$
By Proposition \ref{centralizerGroup},
$$\End_{\Gal(K_0)}(V_{\ell}(A))=\End_{\g_{\ell,A}}(V_{\ell}(A)).$$
This implies that
$$\End_{K_0}(A)\otimes\Q_{\ell}=\End_{\g_{\ell,A}}(V_{\ell}(A)).$$
By Proposition \ref{centralizerLie},
$$\End_{\g_{\ell,A}}(V_{\ell}(A))=\End_{\g_{\ell,A}^0}(V_{\ell}(A)).$$
This implies that
$$\End_{E_0}(A)\otimes\Q_{\ell}=\End_{K_0}(A)\otimes\Q_{\ell}=\End_{\g_{\ell,A}^0}(V_{\ell}(A)).$$
By Remark \ref{Liecentral} applied to $K_0$ and $E_0$ (instead of
$K$ and $E$),
$$\End_{\g_{\ell,A}^0}(V_{\ell}(A))\supset
\End_{\Gal(E_0)}(V_{\ell}(A)).$$ So, we get
$$\End_{E_0}(A)\otimes\Q_{\ell}=\End_{\g_{\ell,A}^0}(V_{\ell}(A))\supset
\End_{\Gal(E_0)}(V_{\ell}(A))\supset
\End_{E_0}(A)\otimes\Q_{\ell},$$ which implies that
$$\End_{E_0}(A)\otimes\Q_{\ell}=\End_{\g_{\ell,A}^0}(V_{\ell}(A))=
\End_{\Gal(E_0)}(V_{\ell}(A))= \End_{E_0}(A)\otimes\Q_{\ell}.$$ In
particular,
$$\End_{E_0}(A)\otimes\Q_{\ell}=\End_{\Gal(E_0)}(V_{\ell}(A)).$$ Now
Lemma \ref{prod}.1 applied to $X=Y=A$ and to $E_0$ (instead of $K$)
implies that
$$\End_{E_0}(A)\otimes\Z_{\ell}=\End_{\Gal(E_0)}(T_{\ell}(A)).$$ It
follows from Remark \ref{extn} (applied to $E_0/E$ instead of $L/K$)
that $$\End_{E}(A)\otimes\Z_{\ell}=\End_{\Gal(E)}(T_{\ell}(A)).$$
Again, Lemma \ref{prod}.1 tells us that
 $$\End_{E}(A)\otimes\Q_{\ell}=\End_{\Gal(E)}(V_{\ell}(A)).$$

\end{proof}

\begin{lem}[of Clifford]
\label{clifford} Let $G$ be a group and $H$ its normal subgroup. Let
$W$ be a vector space of finite positive dimension over a field $k$.
Let $\rho: G \to \Aut_k(W)$ be a semisimple (completely reducible)
linear representation of $G$. Then the corresponding $H$-module $W$
is also semisimple.
\end{lem}

\begin{proof}
Let us split $W$ into a finite direct sum $W=\oplus W_i$ of simple
$G$-modules $W_i$. By Theorem (49.2) of \cite{CR}, the corresponding
$H$-modules $W_i$ are semisimple. This implies that the $H$-module
$W$ is a direct sum of semisimple $H$-modules $W_i$'s and therefore
is also semisimple.
\end{proof}

\begin{prop}
\label {semisimplicity}
 Let $L/K$ be a finite or infinite Galois
extension of $K$. If the $\Gal(K)$-module $V_{\ell}(A)$ is
semisimple then the $\Gal(L)$-module $V_{\ell}(A)$ is semisimple.
\end{prop}

\begin{proof}
 Since $L/K$ is Galois, the subgroup $\Gal(L)$ of $\Gal(K)$ is normal. Now
 the result follows from Lemma \ref{clifford}.
\end{proof}

\section{Homotheties, centralizers and semisimplicity}
\label{central}

\begin{thm}[of Bogomolov]
\label{bogom} Suppose that $K$ is a field that is finitely generated over  $\Q$
and $\ell$ is a prime. Let $A$ be an abelian variety of positive dimension over
$K$, Then $\g_{\ell,A}$ contains the homotheties $\Q_{\ell}\I$.
\end{thm}

\begin{proof}
When $K$ is a number field, this assertion was proven by Bogomolov
\cite{Bogomolov1,Bogomolov2}. The case of arbitrary finitely
generated $K$ is also known \cite[Sect. 1, p. 2]{SerreLR} and
follows from the number field case with the help of Serre's variant
of the Hilbert irreducibility theorem for infinite extensions
(\cite[Sect. 1]{SerreLR}, \cite[Sect. 10.6] {SerreMW}, \cite[Prop.
1.3 on pp. 163--164]{Noot}). Indeed, there exist a number field $F$,
an abelian variety $B$ over $F$ with $\dim(A)=\dim(B)$ and an
isomorphism of $\Z_{\ell}$-modules $u:T_{\ell}(A)\cong T_{\ell}(B)$
such that $u^{-1} G_{\ell,B,F}u=G_{\ell,A,K}$. Extending $u$ by
$\Q_{\ell}$-linearity, we obtain the isomorphism of
$\Q_{\ell}$-vector spaces $V_{\ell}(A) \cong V_{\ell}(B)$, which we
still denote by $u$. Clearly,
$$u^{-1} \g_{\ell,B}u=\g_{\ell,A}.$$
Since $F$ is a number field, $\g_{\ell,B}$ contains all the
homotheties, which implies that $\g_{\ell,A}$ also contains all the
homotheties. This ends the proof.
\end{proof}

{\sl Proof of Theorem \ref{main}.} In light of Remark \ref{enlarge},
we may and will assume that $X$ and $Y$ are defined over $K$. Let us
put $A=X\times Y$. Since $K$ is finitely generated over $\Q$, every
finite algebraic extension $K^{\prime}$ of $K$ is also finitely
generated over $\Q$.
 By the theorem of Faltings \cite{F1,F2}, the injection
$$\tilde{i}_{A,A,K^{\prime}}:\End_{K^{\prime}}(A)\otimes\Q_{\ell} \twoheadrightarrow
\End_{\Gal(K^{\prime})}(V_{\ell}(A))$$ is bijective.  Thanks to
Theorem \ref{bogom}, we know that $\g_{\ell,A}$ contains the
homotheties $\Q_{\ell}\I$. Now the desired result follows from
Theorem \ref{main0} combined with Lemma \ref{prod}.

\begin{cor}
\label{maincor} Suppose that $K$ is a field that is finitely
generated over  $\Q$ and $\ell$ is a prime. Let  $E_1$ be a field
extension of $K$ that lies in $K(\ell)$. If $X$ and $Y$ are abelian
varieties over $E_1$ then the natural embedding of
$\Z_{\ell}$-modules
$$\Hom_{E_1}(X,Y)\otimes \Z_{\ell} \hookrightarrow \Hom_{\Gal(E_1)}(T_{\ell}(X),
T_{\ell}(Y))$$ is bijective.
\end{cor}

\begin{proof}
Since $K(\ell)/E_1$ is a Galois extension, the result follows from
Theorem \ref{main} combined with Remark \ref{extn} applied to
$K(\ell)/E_1$.

\end{proof}

\begin{thm}
Suppose that $K$ is a field that is finitely generated over  $\Q$
and $\ell$ is a prime. Let $L/K$ be a finite or infinite Galois
extension. (E.g., $L=K(\ell)$.)  Let $A$ be an abelian variety of
positive dimension over $K$.  Then the $\Gal(L)$-module
$V_{\ell}(A)$ is semisimple.
\end{thm}

\begin{proof}
Faltings \cite{F1,F2} proved that the $\Gal(K)$-module $V_{\ell}(A)$
is semisimple. This also covers the case when $L/K$ is a finite
(Galois) extension. The case of infinite $L/K$ follows from
Faltings' result combined with Proposition \ref{semisimplicity}.

\end{proof}

\begin{thm}
Suppose that $K$ is a field that is finitely generated over  $\Q$ .
Let $L/K$ be a finite or infinite Galois extension.  Let $A$ be an
abelian variety of positive dimension over $K$.  Then the
$\Gal(L)$-module $A_{\ell}$ is semisimple for all but finitely many
primes $\ell$.
\end{thm}

\begin{proof}
For all but finitely many primes $\ell$ the $\Gal(K)$-module
$A_{\ell}$ is semisimple. Indeed, when $K$ is a number field, this
assertion is contained in Corollary 5.4.3 on p. 316 of
\cite{ZarhinInv} (the proof is based on results of Faltings
\cite{F1}). The same proof works over arbitrary fields that are
finitely generated over $\Q$, provided one replaces the reference to
Prop. 3.1 of \cite{F1} by the reference to the corollary on p. 211
of \cite{F2}. Since $L/K$ is Galois, $\Gal(L)$ is a normal subgroup
of $\Gal(K)$. Now the desired result follows from Lemma
\ref{clifford} applied to $k=\F_{\ell}, W=A_{\ell}$ and $G=\Gal(K),
H=\Gal(L)$.
\end{proof}

\end{document}